\documentclass[11pt]{article}
\usepackage[english]{babel}
\usepackage[latin1]{inputenc}
\usepackage{latexsym,amsfonts,amsmath,amssymb,euscript,epsfig,graphicx,color,float}
\usepackage[active]{srcltx}
\usepackage[T1]{fontenc}
\title{A differential game on Wasserstein space. Application to weak approachability with partial monitoring\thanks{This research was partially by research contract AFOSR-FA9550-18-1-0254.}}

\author{Vianney Perchet \thanks{CMLA, ENS Paris-Saclay \& Criteo Research, France. (vianney.perchet@normalesup.org). } \and Marc Quincampoix \thanks{Laboratoire de Math\'ematiques de Bretagne Atlantique (CNRS UMR 6205), 6, Avenue Victor Le Gorgeu, 29200 Brest, France. (Marc.Quincampoix@univ-brest.fr ). }}

\usepackage{amsfonts}

\newcommand\cY{\mathcal{Y}}
\newcommand\cX{\mathcal{X}}
\newcommand\cA{\mathcal{A}}
\newcommand\cZ{\mathcal{Z}}

\newcommand\DD{\Delta_\delta}
\newcommand\cH{\mathcal{H}}

\newcommand\bx{\mathbf{x}}
\newcommand\by{\mathbf{y}}
\newcommand\bz{\mathbf{z}}

\newcommand\cG{\mathcal{G}}
\newcommand\cK{\mathcal{K}}

\newcommand\cS{\mathcal{S}}

\newcommand\cB{\mathcal{B}}
\newcommand\N{\mathbb{N}}
\newcommand\R{\mathbb{R}}

\renewcommand{\d}{\mathrm{d}\,}
\newcommand{\Lip}{\mathrm{Lip}}

 \newcounter{hypothese}
\newtheorem{hypo}[hypothese]{Assumption}
 \newcounter{theo}
\newtheorem{definition}[theo]{Definition}
\newtheorem{proposition}[theo]{Proposition}
\newtheorem{lemma}[theo]{Lemma}
\newtheorem{theorem}[theo]{Theorem}
\newtheorem{corollary}[theo]{Corollary}

\begin{document}
\maketitle

\begin{abstract}
Studying continuous time counterpart of some discrete time dynamics is now a standard and fruitful technique, as some properties hold in both setups. In game theory, this is usually done by considering  differential games on  Euclidean spaces. This allows to infer properties on the convergence of values of a  repeated game, to deal with the various concepts of approachability, etc. In this paper, we introduce a specific but quite abstract differential game defined on the Wasserstein space of probability distributions and we prove the existence of its value. Going back to the discrete time dynamics, we derive results on weak approachability with partial monitoring: we prove that any set satisfying a suitable compatibility condition is either weakly approachable or weakly excludable.
We also obtain that 
the  value for differential games with nonanticipative strategies is the same that those defined with a new concept of strategies very suitable to make links with repeated games.

\end{abstract}

\section{Introduction}

Blackwell's approachability \cite{Bla56} is a core concept in repeated games \cite{AumMas55,Tomala,AbeBarHaz11}. It is defined in two-player repeated games where the stage outcome is a vector in $\R^d$, possibly representing $d$ different criteria to optimize simultaneously. Both players aim at controlling the time average outcome.  The objective of player 1 is that the time average vectorial payoff converges to some fixed target set $E \subset \R^d$. If he can ensure that objective, then the target set is called approachable. The objective of player 2 is to prevent this convergence.

The motivation behind approachability theory is twofold as it can be applied both in game theory and in machine learning (more specially in online learning). Indeed, it is  a standard tool in game theory, as it can be used  to construct optimal strategies in repeated games with incomplete information \cite{AumMas55} or to construct equilibria in multi-player repeated game \cite{Tomala}. In machine learning,  it offers a clean and elegant solution to online multi-criteria optimization problems.  More precisely, one of the most important class of problems called regret minimization \cite{Bla56a,AbeBarHaz11}, as well as other online learning criteria such as calibration \cite{Daw85,FosVoh97,Per14,Per15,MPS14} are special cases of approachability. We refer to \cite{CBL,Per14} for surveys and textbooks on online learning and connections between approachability and other concepts.

\subsection{Weak approachability in Euclidean space, with full monitoring}\label{SE:WeakFull}
We assume that the action sets of  player 1 and 2 are convex and compact subset of  $\R^a$ and $\R^{b}$ denoted respectively by $X \subset \R^a$  and $Y \subset \R^{b}$. The outcome is  defined trough a bi-linear mapping $g(x,y)=xAy:=\sum_{i, j} x[i]y[j]A_{i,j} \in \R^d$ where $A_{i,j} \in \R^d$ and $x=(x[i])_{ i \in \{1,\ldots,a\}}$, $y=(y[j])_{ j \in \{1,\ldots,b\}}$. We denote by $\cG=\{g(x,y), x \in X, y \in Y\}$ the range of $g$ and by $g_m:=g(x_m,y_m)$ the outcome at stage $m \in \N$ generated by the choices of $x_m \in X$ and $y_m \in Y$.

We assume that the length of the game is finite, known  and equal to  $n \in \N$. A strategy of  player 1 is then a mapping from $\bigcup_{k=0}^{n-1}\left(X \times Y\right)^k$ into $X$ (and into $Y$ for  player 2).

\medskip

In this vectorial framework,  objectives are represented by some exogenous closed  set $E \subset \R^d$. Player 1 aims at making the average outcome $\overline{g}_n:=\frac{1}{n}\sum_{m=1}^n g_m$ converge to $E$ and player 2 aims at preventing it. Stated otherwise, player 1 aims at minimizing the distance $d(\overline{g}_n,E)$ of $\overline{g}_n$ to $E$, where $d(z,E)=\inf_{\omega\in E} \|z-\omega\|$; conversely, player 2 aims at maximizing this distance.

\medskip
\begin{definition}\label{DEF:WA}
A closed set $E \subset \R^d$ is weakly approachable by the first player if for every $\varepsilon >0$, there exists $N \in \N$, such that for every $n \geq N$, there exists a strategy of the first player $\sigma_n$  such that, no matter the strategy $\tau$ of the second player,  $d(\overline{g}_n,E) \leq \varepsilon$.

\medskip

A closed set $E \subset \R^d$ is weakly excludable if player 2 can weakly approach the complement of some $\eta$-neighborhood of $E$ (with $\eta>0$ small enough).
\end{definition}

Vieille \cite{Vie92}  proved the following  conjecture of Blackwell \cite{Bla56}

\medskip

\begin{proposition}\label{PR:Vieille}
Every closed set $E \subset \R^d$ is either weakly approachable or weakly excludable.
\end{proposition}

\medskip

Informally speaking, this result says that the value\footnote{Because of the lack of convexity of the payoff mapping, the value is only guaranteed to exist in the equivalent framework where players can choose actions at random.}  of the zero-sum game with payoff $d(\overline{g}_n,E)$ converges as the horizon $n \in \N$ increases, see, e.g., \cite{Per14}. The main insight behind the proof of Vieille \cite{Vie92}  consists  in seeing the $n$-stage game, for $n$ large enough, as a close approximation of a specific zero-sum differential game. As a consequence, the values of the repeated games  converge to the value of the differential game. Recently,   connections between differential games and repeated games have been fruitfully exhibited  \cite{CarLarSor12,AssQuiSor09,Per14}.

\subsection{A Basic introduction of differential games}\label{basic}

Differential games were introduced in \cite{I, Pont}. Here we consider two-player zero-sum differential games which dynamics  are \begin{equation}\label{equ'}  g'(s) = f(s, g(s), x(s), y(s)) , \; s \in  [0,1] \end{equation}  and which payoff has the form   \begin{equation*}  {\cal J } := \ell ( g(1)). \end{equation*}  The state variable is $ g \in \R ^d$. The players acts on the system by choosing measurable controls : Player 1 wants to minimize the payoff  ${\cal J }$ by choosing the control $ x : [0,1] \mapsto X $ while Player 2 tries to maximize ${\cal J }$ by choosing the control $ y : [0,1] \mapsto Y $ ($X$ and $Y$ are given compact subsets of some finite dimensionnal space).
 The functions $ f : [0,1] \times \R ^d\times   X \times Y  \mapsto \R ^ d$ and $ \ell : \R ^d \mapsto \R$  are supposed to be Lipschitz continuous and bounded.

In view of his objective (minimization or maximization) of the payoff, each player chooses his control knowing the past actions of his opponents. This is precisely expressed by the notion of strategies we explain now.

Let us denote by  $\cX$ the set of measurable controls $ x : [0,1] \mapsto X $ of player 1. Similarly  $\cY(s_0)$ denotes the set of player 2 controls.  A nonanticipative strategy for player 1 is a map $ \alpha : \cY \mapsto \cX $ which associates to any a control $y$ chosen by player 2 a control $ x$ of player 1 in a nonanticipative way i.e : For any $t \in [0,1] $ if two controls $y_1$ and $y_2$ coincide almost everywhere on $[0,t]$ then $ \alpha (y_1)$ and $ \alpha (y_2)$ also  coincide almost everywhere on $[0,t]$. An nonanticipative strategy $ \beta  : \cX \mapsto \cY $  for player 2 is similarly defined.

 For $(s_0, g_0) \in [0,1] \times \R ^d$ , we denote by $ g ^{s_0, g_0, x, y}(\cdot)  $ the unique solution of (\ref{equ'}) with the initial condition $g(s_0) =g_0$. We define then the following value functions which are the results of the optimal actions of the players
 \begin{eqnarray*}
V^+ (s_0,g_0) := \inf _ \alpha \sup _{y \in \cY}  \ell\Big(g ^{s_0, g_0, \alpha (y) , y}(1)\Big) \\
V^- (s_0,g_0) := \sup _ \beta \inf _{x \in \cX}  \ell\Big(g ^{s_0, g_0, x , \beta (x) }(1)\Big).
\end{eqnarray*}
An important problem concerns the existence of a value of the game namely the validity of the equality $V ^+ = V ^-$. This result has been for instance obtained in \cite{ES} by proving that $V ^+ $ and  $ V ^- $ are Lipschitz continuous and they are both viscosity solution of a partial differential equation (called the Hamilton Jacobi Isaacs equation) which has the uniqueness property. This result is valid under suitable Isaacs' condition (cf \eqref{I} later on).  For more general existence of value results we refer the reader to \cite{BeCQ,CQSP,PQ} (see also \cite{BCQ}).

It could be surprising at the first glance that the above values $V^+ $ and  $V ^-$ are not defined in a symmetric way. One can show that the game can be written into a a normal form when the cost $ {\cal J} $ is regular enough with a little different notion of strategies  (cf Definition \ref{NAD}). Section \ref{SE:Appendix} is devoted to some basic facts for differential games and also contains a new result concerning a new class of strategies.

\subsection{From partial monitoring in Euclidean space to full monitoring in Wasserstein space}
A crucial and implicit assumption in the  model  of Section \ref{SE:WeakFull} is the fact  that players observe at each stage the action chosen by their opponent. This framework is usually called  ``\textsl{with full monitoring}''.

\bigskip

As it is now standard in game theory \cite{AumMas55,Koh75,MerSorZam94,Per11a} and machine learning \cite{Rus99,LugManSto08,ManPerSto11,Per11b,KP} we may assume that player~1 does not necessarily observe the action of player 2 but  only receives signals about it. This  framework is  called  ``\textsl{with partial monitoring}''.   Formally, if actions taken at stage $m$ are $x_m \in X$ and $y_m \in Y$ then the (unknown) stage outcome is $g_m=x_mAy_m$ and the signal observed is $\mu_m:= Sy_m \in \R^k$ where $S$ is a $1\times b$ matrix with  components in $\R^k$.

A strategy of player 1 is then a mapping from $\bigcup_{m=0}^{n-1}\left(X \times \mathcal{S}\right)^m$ into $X$, where  $\cS$ is the range of $S$. On the other hand, we can assume that player 2 has still a full monitoring on player 1's actions, and his strategies are mappings from $\bigcup_{m=0}^{n-1}\left(X \times Y \times \mathcal{S}\right)^m$ into $Y$. As with full monitoring, a closed set $E \subset \R^d$ is weakly approachable if player 1 can ensure that the average payoff is $\varepsilon$-close to $E$, if the length of the game is big enough, see Definition \ref{DEF:WA}.

\medskip

Perchet and Quincampoix \cite{PerQui14} have developed an abstract setup to treat any game with partial monitoring as a game with full monitoring and  outcomes in the Wasserstein space of probability distribution on $X\times \mathcal{S}$. Outcomes are therefore probability measures to be interpreted as the maximal information available to player 1.

 The basic idea relies on the following multi-valued mapping $\mathbf{p}$ defined on $X \times \cS$.
\[ \mathbf{p}(x,\mu)= \Big\{ xAy \ ; \ y  \in Y \ \mathrm{such\ that}\ Sy=\mu \Big\}.
\]
The sequence $(y_1,\ldots,y_n)$  of  actions of player 2 generates a sequence of signals $(\mu_1,\ldots,\mu_n)$. So the only information available to the player is that $g_m$, the outcome at stage $m$, belongs to $\mathbf{p}_m=\mathbf{p}(x_m,\mu_m)$. As a consequence,   to ensure that $\overline{g}_n$ belongs to some set $E$,  the Minkowski average set
$$\overline{\mathbf{p}}_n = \frac{1}{n}\sum_{m=1}^n\mathbf{p}_m := \Big\{\frac{1}{n}\sum_{m=1}^n g_m \, ; \, g_m \in \mathbf{p}_m  \Big\}$$
 must be  included in $E$. Moreover, this inclusion has the following interesting interpretation:
\begin{eqnarray*} \overline{\mathbf{p}}_n  \subset E &\Leftrightarrow \frac{1}{n} \sum_{m=1}^n \mathbf{p}(x_m,\mu_m) \subset E  \Leftrightarrow \frac{1}{n} \sum_{m=1}^n \mathbb{E}_{\delta_{x_m}\otimes\delta_{\mu_m}}[\mathbf{p}] \subset E \Leftrightarrow   \mathbb{E}_{\overline{(x\otimes \mu)}_n}[\mathbf{p}]  \subset E  \\ &\Leftrightarrow \overline{(x\otimes \mu)}_n \in \widetilde{E}:=\Big\{ q \in \Delta(X\times \mathcal{S}) \ \mathrm{s.t} \  \mathbb{E}_{q}[\mathbf{p}]  \in E \Big\},
\end{eqnarray*}
where $ \otimes$ stands for the tensor product, $\delta_x$ is the Dirac mass on $x$ and we introduce the notation $(x\otimes \mu)_m = \delta_{x_m}\otimes\delta_{\mu_m}$. The set $\widetilde{E} \subset \Delta(X\times \mathcal{S})$ we introduced somehow corresponds to the set of ``compatible informations'' with the objectives of players 1, i.e., those that guarantee  the average payoff belongs to the target set $E$. Stated otherwise, the problem of weak-approachability  with partial monitoring of a closed set $E \subset \R^d$ can be rewritten as a problem of weak-approachability of the set $\widetilde{E} \subset \Delta(X\times \cS)$  in the Wasserstein space . This reduction has the following interesting upside: instead of trying to control average of sets in Euclidean space, a player can aim at controlling averages of points, even if they belong to some lifted,  more complex space. Indeed, averages of sets are difficult objects to handle; for instance, even intuitively, it is not clear how to make them ``converge'' to a target set. On the contrary, it is rather intuitive for averages of points (even in a lifted space): one just need to find the next point in the ``direction'' of the target set. Here, the difficulty is to define properly the concept of \textsl{direction}, yet once this is done, controls are easier to construct.

\bigskip

As a consequence, we aim at generalizing the traditional techniques in Euclidean space by introducing some  zero-sum differential games in Wasserstein space. Then we will first  obtain conditions ensuring the existence of a value and then show that the $n$-stage repeated games are close to some limit differential game. This will enable us   to prove that any closed set $\widetilde{E}\subset \Delta(X\times \cS)$ is either weakly-approachable or weakly-excludable.

Notice that this does not imply that any set $E \subset \R^d$ is either weakly approachable or weakly excludable (which is an incorrect statement, see Perchet \cite{Per11a}), but it implies that any set of the form
\[ \{\mathbb{E}_{q}[\mathbf{p}] , q \in \widetilde{E} \} \subset \R^d
\]
is either weakly-approachable or weakly-excludable with partial monitoring.

\subsection{Weak approachability in Wasserstein Space}
We now define formally what we meant in the last section by weak approachability in the Wasserstein space. Let $X \subset \R^a$ and $Z\subset \R^{b}$ be two convex compacts sets of some Euclidian spaces and $\widetilde{E} \subset \Delta(X\times Z)$ be a closed subset of  $\Delta(X \times Z) $, the set of probability distributions over $X\times Z$. This set is equipped  with the Wasserstein quadratic distance $W_2$, which  definition is recalled  in the preliminaries section \ref{prelimi}.

\medskip

The game in discrete time is described as follows. At stage $m\in \N$,  players  choose respectively $\bx_m \in \Delta(X)$ and $\bz_m \in \Delta(Z)$ and these choices induce the stage outcome $\theta_m=\bx_m \otimes \bz_m \in \Delta(X \times Z)$. This outcome $\theta_m$ is observed by both players, i.e., the game is with full monitoring, before stage $m+1$ begins.

\medskip

\begin{definition}
A closed set $\widetilde{E} \subset \Delta(X\times Z)$ is weakly approachable by  player 1 if for every $\varepsilon >0$, there exists $N \in \N$, such that for every $n \geq N$, he has a strategy $\sigma_n$  such that, no matter the strategy $\tau$ of  player 2,
 $$W_2(\overline{\theta}_n,\widetilde{E}):= \inf_{q \in \widetilde{E}} W_2(\overline{\theta}_n,q) \leq \varepsilon.$$

\medskip

A closed set $\widetilde{E} \subset \Delta(X\times Z)$ is weakly excludable by player 2 if he can weakly approach the complement of some $\eta$-neighborhood of $\widetilde{E}$, with $\eta>0$.
\end{definition}

\subsection{Organization of the paper and main results}
The remaining of the paper is divided in two main parts. Section \ref{SE:Euclidean} is devoted to the study of differential and repeated games in Euclidean space to get some intuitions, and we start working  in the  Wasserstein space of probability measures  in  Section \ref{SE:Wasserstein}, after a preliminary section on Wasserstein distance.

In Section \ref{SE:Euclidean}, we basically recover the main result of Vieille but with an alternative proof and with new concepts of strategies that we purposely introduced. One might wonder why we bother proving again such an elegant result with a longer and maybe more intricate proof. The first reason is that the Euclidean framework is obviously more natural and more intuitive than the Wasserstein space of probability measures. But more importantly,  all the proofs we give in the former setup can be  generalized at no cost to the later. Unfortunately, it was not the case of the techniques of Vieille \cite{Vie92}, based notably on results of Flemming \cite{Fle64a,Fle64b} or some of the ideas appearing in some other differential games with dynamics in Wasserstein space \cite{CarQui08, JQ, MQ}.

For the sake of clarity, we therefore chose to decompose the main arguments into those that hold no matter the ambiant space (i.e., either in Euclidean or Wasserstein spaces, the proof being stated in the former as it is  more intuitive) and those that are true in the Wasserstein space; they are described in Section \ref{SE:Wasserstein}.

\section{Weak approachability with full monitoring through a differential game} ~

\subsection{A Differential Game with Non Anticipative strategies with delay}

\label{SE:Euclidean} As it becomes more and more popular in repeated game theory, we represent the $n$-stage repeated game as a discretization (or an approximation) of some differential game. \medskip

  Given the fixed horizon $n \in \N$ and $m \leq n$, the following equation describes the evolution of average payoffs in discrete time.
\[ \overline{g}_{m+1}^n =\overline{g}_{m}^n + \frac{1}{m+1}(x_{m+1}^nAy_{m+1}^n-\overline{g}_{m}^n).
\]
The continuous analogue of the above discrete equation is  the  following  differential equation:
\begin{equation}\label{EQ:DynDiffDiscreteFull} \dot g(s) =\frac{1}{s}\Big( \bx(s)A\by(s)-g(s)\Big), \quad \forall s \in [s_0,1]\ \mathrm{\ and\ } g(s_0)=g_0 \in \cG  \subset \R^d
\end{equation}
for some $s_0 >0$. Its  solution is given by the following integral equation
 \begin{equation}\label{EQ:g}g(s)=\frac{s_0}{s}g_0+\frac{1}{s} \int_{s_0}^s\bx(t)A\by(t)dt.\end{equation}

  A control $\bx$  of the player 1  is a measurable map from $[s_0,1]$ to $X$; the set of such controls  is denoted by $\cX(s_0)$ (and, similarly, $\cY(s_0)$ for player 2).
For any $ s_0 \in (0, 1] $, $ g_ 0 \in \R ^d $ and
  $(\bx, \by) \in \cX (s_0) \times \cY (s_0) $  we denote by $ s \mapsto g^{s_0,g_0,\bx, \by}(s) $ the unique solution to (\ref{EQ:g}).

Let us recall the notion of Nonanticipative Strategies with Delays (in short: NAD Strategies). \begin{definition} \label{NAD} Given $s _0 \in \R$, a NAD for player 1 is a map $ \alpha : \cY(s_0) \mapsto \cX(s_0) $ such that there exists a subdivision $ t_0:= s_0 <t_1 < \ldots t_N :=1 $ of the interval $ [s_0 ,1] $ such that for any $ k =0, 1 \ldots N-1 $ if $ \by _1 ( \cdot)   $ and $ \by _2 ( \cdot) $ coincide almost surely on $ [s_0 , t_k]$ then the controls  $ \alpha (\by _1 ) ( \cdot)   $ and $ \alpha (\by _2)  ( \cdot) $ coincide almost surely on $ [s_0 , t_{k+1}]$. The set of such nonanticipative strategies  $ \alpha $ for player 1 is denoted by $ \mathcal{A}(s_0)$. We define in a similar way the set $ \mathcal{B}(s_0) $ of nonanticipative strategies $ \beta $ for  player 2.\end{definition}

One interest of such strategies lies on the fact that one can associate a trajectory to a pair of strategies due to the following result, cf, e.g., \cite{CarQui08}.
\begin{lemma} \label{rep} For any pair $(\alpha, \beta) \in  \mathcal{A} (s_0) \times  \mathcal{B}   (s_0 ) $ there exists a unique pair of control   $(\bx, \by) \in \cX (s_0) \times \cY (s_0) $  such that
$$ \alpha ( \bx) = \by  \mbox{ and } \beta ( \bx ) = \by .$$

\end{lemma}
So for  $(\alpha, \beta) \in  \mathcal{A} (s_0) \times  \mathcal{B}   (s_0) $
 we define   $  g^{s_0,g_0,\alpha, \beta } = g^{s_0,g_0,\bx, \by} $ where  $(\bx, \by) $ is associated with  $(\alpha, \beta) $ by Lemma \ref{rep}.

Coming back to the dynamics (\ref{EQ:DynDiffDiscreteFull}) and considering the specific loss $\ell : \R^d \to \R_+$ defined by $\ell(z) : =d(z,E)$, the distance to the closed set $E$, we now define the value functions of the game.

The upper-value of the differential game is given by
\[V^+(s_0,g_{0}):=\inf_{\alpha \in  \mathcal{A} ( s_0  ) }\sup_{\beta \in  \mathcal{B} ( s_0  ) }\ell\Big(g^{s_0,g_0,\alpha,\beta}(1)\Big).\]
while the  the lower-value is
\[V^-(s_0,g_{0}):=\sup_{\beta \in  \mathcal{B} ( s_0  ) }\inf_{\alpha \in  \mathcal{A}( s_0 ) }\ell\Big(g^{s_0,g_0,\alpha,\beta}(1)\Big).\]
for every $s_0 \in (0,1]$ and $g_0 \in \R ^d$.

Observe  also that Lemma \ref{rep} yields \[V^+(s_0,g_{0})=\inf_{\alpha \in  \mathcal{A} ( s_0  )}\sup_{\by}\ell\Big(g^{s_0,g_0,\alpha(\by),\by}(1)\Big)\ \mbox{ and } V^-(s_0,g_{0})=\sup_{\beta \in  \mathcal{B} ( s_0 )}\inf_{\bx}\ell\Big(g^{s_0,g_0,\bx,\beta(\bx)}(1)\Big) \] which is the same definition of values that those given in subsection  \ref{basic}. Some classical results on such differential games are recalled in the Appendix.

Because of the dynamics we consider, the game has a value in NAD strategies :
\[ \forall \, s_0 >0, \forall g _0 \in \R^d , \;   V^-(s_0,g_0)  = V^+(s_0,g_0)  .\] Moreover the common value -denoted by $V$-  is the unique Lipschitz continuous viscosity solution of
\begin{equation}\label{HJB} \left\{\begin{array}{ll} \frac{\partial V}{\partial s}(s,g)+H\left(s,g,\frac{\partial V}{\partial g}(s,g)\right) =0 & \mathrm{for\ all }\ (s,g) \in [s_0,1]\times \R^d\\
V(1,g)=\ell(g) & \mathrm{for\ all }\ g \in \R^d\end{array}\right.
\end{equation} where
\[H(s,g,p):=\sup_{y\in Y} \inf_{x \in X}p.\frac{1}{s}\Big\{xAy-g\Big\}= \inf_{x \in X}\sup_{y\in Y}p.\frac{1}{s}\Big\{xAy-g\Big\}\]

In view of a deeper analysis between the above discrete game and the differential game introduced, we are led to introduce a smaller class of strategies in the following Section. We mention here that the  results of the following Section \ref{SE:NADCStrat}, stated for this new class of strategies, also hold for the classical NAD strategies. Yet  this new concept of strategy is, first,  conceptually simpler and, second,   can be directly connected to a strategy in a repeated game (see Section \ref{SE:FromDiffToRep}).

\subsection{A new concept of strategies adapted to the discrete-continuous time approximation }\label{SE:NADCStrat} In this section, we develop the new concept of  strategies, more adapted to the discrete-continuous time approximation than the  existing ones \cite{Fle64a,Fle64b}.

\bigskip

 \medskip

\begin{definition}\label{DF:NADCstrat}
A non-anticipative with delay piecewise-constant strategy of player 1 is a mapping $\alpha$ from $\cY(s_0)$ to $\cX(s_0)$  satisfying the following properties:
 \begin{itemize}
 \item[1)] There exists some integer $N \in \N$ such that, for any control $\mathbf{y} \in \cY(s_0)$, $\alpha(\mathbf{y})$ is constant on $[m/N,(m+1)/N]$, for all $m \in \{0, \ldots, N-1\}$
 \item[2)] The strategy is non-anticipative with delay: if $\mathbf{y}(s)=\mathbf{y}'(s)$ for all $s\in [s_0, m/N]$, then $\alpha(\mathbf{y})(\cdot)=\alpha(\mathbf{y}')(\cdot)$ on $[s_0,(m+1)/N]$
 \end{itemize}
 We denote by $ \mathcal{A}_{NADC}(s_0)$ the set of such strategies of the player 1 and, similarly,  those of player 2 by $ \mathcal{B}_{NADC}(s_0)$.
\end{definition}

\medskip

Since every strategy $\alpha \in \cA_{NADC}(s_0)$ and $\beta \in \cB_{NADC}(s_0)$ are non-anticipative with delay,  there exists a unique pair of control $\bx \in \cX(s_0)$ and $\by \in \cY(s_0)$ such that $\alpha(\by)=\bx$ and $\beta(\bx)=\by$.  We can define as usual the game in normal form and the values.

\medskip

\begin{definition}
The upper-value of the game is defined for every $s_0 \in (0,1]$ and $g_0 \in \cG$ by
\[V_{NADC}^+(s_0,g_{0})=\inf_{\alpha \in \cA_{NADC}(s_0) }\sup_{\beta  \in \cB_{NADC}(s_0)}\ell(g^{s_0,g_0,\alpha,\beta}(1))=\inf_{\alpha \in \cA_{NADC}(s_0) }\sup_{\by}\ell(g^{s_0,g_0,\alpha(\by),\by}(1)).\]
Similarly, the lower-value is defined for every $s_0 \in (0,1]$ and $g_0 \in \cG$  by
\[V^-_{NADC}(s_0,g_{0})=\sup_{\beta\in \cB_{NADC}(s_0) }\inf_{\alpha \in \cA_{NADC}(s_0) }\ell(g^{s_0,g_0,\alpha,\beta}(1))=\sup_{\beta\in \cB_{NADC}(s_0) }\inf_{\bx}\ell(g^{s_0,g_0,\bx,\beta(\bx)}(1)).\]
\end{definition}
\medskip

It always holds that $V_{NADC}^- \leq V_{NADC}^+$ and, because of the dynamics we consider, these mappings are  regular:
\medskip

\begin{lemma}\label{LM:ValuesLipschitz} Let $\kappa>0$ be  a uniform bound on $g_0$ and  $\|xAy\|$, then
\begin{itemize}
\item[1)] $V_{NADC}^+(\cdot,g_0)$ and $V_{NADC}^-(\cdot,g_0)$  are  $2\kappa$-Lipschitz and can be extended to $[0,1]$.
\item[2)] $V_{NADC}^+(0,g_0)$ and $V_{NADC}^-(0,g_0)$ are  independent of $g_0$
\item[3)] $V_{NADC}^+(s,\cdot)$ and $V_{NADC}^-(s,\cdot)$ are $2\kappa s$-Lipschitz
\end{itemize}
\end{lemma}

\medskip

\textbf{Proof.}
Since for every control  $\mathbf{x} \in \cX(s_0)$ and $\by \in \cY(s_0)$, and any $s_1\geq s_0$, it holds\begin{eqnarray*} g^{s_0,g_0,\bx,\by}(1)-g^{s_1,g_0,\bx,\by}(1)&=&s_0g_0+\int_{s_0}^1\bx(t)A\by(t) dt-s_1g_0-\int_{s_1}^1\bx(t)A\by(t)\\&=&(s_0-s_1)g_0+\int_{s_0}^{s_1}\bx(t)A\by(t)dt\, ,
\end{eqnarray*}
 one immediately obtains that $V_{NADC}^+(\cdot,g_0)$, $V_{NADC}^-(\cdot,g_0)$  are $2\kappa$-Lipschitz, where $\kappa$ is a uniform bound on $\|g_0\|$ and $\|xAy\|$.  As a consequence, they can be uniquely extended to a $2\kappa$-Lipschitz mapping on $[0,1]$.

 \bigskip
Since $\ell$ is 1-Lipschitz and
  \[\|g^{s_0,g_0,\bx,\by}(1)-g^{s_0,g'_0,\bx,\by}(1)\|=s_0\|g_0-g'_0\|\leq 2\kappa s_0,
 \]
we  obtain that both $V_{NADC}^+(s_0,\cdot)$ and $V_{NADC}^-(s_0,\cdot)$ are $2\kappa s_0$-Lipschitz. As a consequence, the limit when $s_0$ goes to zero of  $V_{NADC}^+(s_0,g_0)$ is independent of $g_0$. $\hfill \Box$

\bigskip


It is worth pointing out that due to the specific form of the dynamics (\ref{EQ:DynDiffDiscreteFull}) we are considering, the above lemma \ref{LM:ValuesLipschitz} is only valid for $s_0 >0$.

\bigskip
It remains to show that  a value exists, i.e., that $V_{NADC}^+(s_0,g_0)=V_{NADC}^-(s_0,g_0)$.  This is due to  results valid on a more general context than the dynamics (\ref{EQ:DynDiffDiscreteFull}). We only sketch the proof, details can be found in the appendix.

\begin{proposition} The game has a value in NADC strategies which coincides with the value of the game in NAD strategies, i.e., for every $s_0 >0$ and $g_0 \in \R^d$,
\[ V^-(s_0,g_0)  = V_{NADC}^-(s_0,g_0)=V_{NADC}^+(s_0,g_0) =V^+(s_0,g_0)=V(s_0,g_0)  .\]
\end{proposition}

\textbf{Proof.}
Because  of Sion \cite{Sio58} minmax theorem, $H^+=H^-$.  Since both  equations (\ref{HJ+}) and (\ref{HJ-}) - stated in the appendix -  reduce to the equation (\ref{HJB}),  Isaacs condition (\ref{I}) holds true. The result is a  direct consequence of Theorem \ref{egalS} of the appendix.
$\hfill \Box$

\medskip

From now on we will shortly  denote by $V$ the value of the game and since  $V(0,\cdot)$ is constant, we simply denote it by $V(0)$.
\subsection{From strategies in differential game to strategies in repeated game}\label{SE:FromDiffToRep} We introduced the new concept of strategies as they are more adapted to the discretization of differential games into repeated games. We now explain this claim through the following lemma.
\medskip

\begin{lemma}\label{LM:FromDiffToRep}
A strategy $\alpha \in \cA_{NADC}(s_0)$, whose delay is $1/N$ (see Definition \ref{DF:NADCstrat}), naturally induces, in the $n$-stage repeated game with $n\geq N$,  a strategy $\sigma_{\alpha,n}$  satisfying \[ \left\|\overline{g}_n - g^{s_0,g^*_0,\alpha(\mathbf{y}),\mathbf{y}}(1)\right\| \leq  \frac{2N\kappa}{n},\]
where $\by$ is the continuous piece-wise constant version of the strategy of player 2 in the $n$ stage repeated game, and $g_0^*$ is some specific point in $\cG$.
\end{lemma}
\medskip

\textbf{Proof.}  Let $1/N$ be the delay of $\alpha$ given by Definition \ref{DF:NADCstrat} and $m^* \in \N_*$ be such that $m^*-1<s_0N\leq m^*$. We denote by $x_0$ the value of $\alpha(\mathbf{y}')(s_0)$ which is independent of the control of player 2, because $\alpha$ is non-anticipative with delay.
We also assume that $n \geq N$ and we let $k \geq 1 $ be such that $n \geq kN+r$ with $r < N$.

\medskip
We  construct the strategy $\sigma_{\alpha,n}$ as follows:
\begin{enumerate}
\item During the first $km^*$ stages, play some arbitrary action $x_0$.
\item The strategy of player 2 in discrete time generates a sequence $y_1,y_2,\ldots,y_{kN} \in Y$. Define a  control $\by$ in continuous time by $\by(\frac{m-1}{kN}+s)=y_m$ for any $s  \in [0,\frac{1}{kN})$
\item At stage $m \in \{km^*+1,\ldots,kN\}$, $\sigma_{\alpha,n}$ dictates to play $\alpha(\by)[\frac{m-1}{kN}]$
\item During the last $r$ stages, play again arbitrarily.
\end{enumerate}
If $n <N$, then $\sigma_{\alpha,n}$ is defined arbitrarily.

\bigskip

By construction of $\sigma_{\alpha,n}$ it immediately reads that $\overline{g}_{kN}=g^{s_0,g^*_0,\alpha(\mathbf{y}),\mathbf{y}}(1)$ where $\mathbf{y}$ is the control defined above and $g^*_0=g^{1/kN,x_0Ay_1,x_0,\mathbf{y}}(s_0)$. As a consequence,
\[ \left\|\overline{g}_n - g^{s_0,g^*_0,\alpha(\mathbf{y}),\mathbf{y}}(1)\right\| \leq  \frac{2r\kappa}{n} \leq \frac{2\kappa}{m^*},
\]
hence the result.$\hfill \Box$

\bigskip

We can finally state and recover  the main result of Vieille.
\medskip

\begin{theorem}\label{TH:Vieille} A closed set $E$ is either weakly approachable or weakly excludable. More precisely, it is weakly approachable by  player 1 if
 $V(0)=0$ and weakly excludable by player 2 if  $V(0)>0$.
\end{theorem}
\medskip

\textbf{Proof.}  Assume that $V(0)=0$, let $\varepsilon >0$ be fixed and  $s_0 = \varepsilon/2\kappa$. Lemma~\ref{LM:ValuesLipschitz} implies that $V(s_0,g_0) \leq \varepsilon$ for any $g_0 \in \cG$. Let  $\alpha \in \cA_{NADC}$ be any $\varepsilon$-optimal strategy in this game and $\sigma_{\alpha,n}$ the associated strategy provided by Lemma \ref{LM:FromDiffToRep}, then

\[ \ell\Big(g^{s_0,g^*_0,\alpha(\by),\by}(1)\Big)\leq \ell\Big(g^{s_0,g_0,\alpha(\by),\by}(1)\Big)+s_0\|g_0-g_0^*\| \leq \left(V(\delta,g_0)+\varepsilon\right)+\varepsilon \leq 3\varepsilon\, .
\]
Assuming that $n=kN+r$ with $2\kappa/k \leq \varepsilon$, we finally obtain that the distance from $\overline{g}_n$ to $E$ is smaller than $4\varepsilon$. Thus, $E$ is weakly approachable.

\bigskip

If $V(0)=\eta>0$, then the same proof gives the fact that the complement of the $\eta$-neighborhood of $E$ is weakly excludable by player 2, so $E$ is weakly excludable. $\hfill \Box$

\section{Weak approachability with partial monitoring through a differential game on Wasserstein space}\label{SE:Wasserstein}

As mentioned in the introduction, we aimed at generalizing the precedent results obtained in a standard Euclidean space to the space of probabilities measures, embedded with the Wasserstein distance. First, we provide some reminder and notations on the Wasserstein distance (and space) and then we describe the associated differential game.

\subsection{Preliminaries on Wasserstein distance }\label{prelimi}

We define in this section  the distance Wasserstein distance  $W_2$ already mentioned in the introduction. We also introduce some material that will be used in the sequel. The reader can refer for this part to the books \cite{AmbGigSav05, Dud89, sant,  Vil03}. For this section only, let us denote by $\cK$ a compact set of some Euclidean space, whose Euclidean  norm is denoted by $\|\cdot\|$.

For every $\mu$ and $\nu$ in $\Delta\left(\cK\right)$, the set  of probability measures on $\cK$, the (squared) Wasserstein distance between $\mu$ and $\nu$ is defined by:

\begin{equation}\label{defwass1}W^2_2(\mu,\nu):=\inf_{\gamma \in \Pi(\mu,\nu)}\int_{ \cK^2}\|x-y\|^2\mathrm d\gamma(x,y)\end{equation}
 where $\Pi(\mu,\nu)$ is the set of probability measures $\gamma \in \Delta\left(\cK\times\cK \right)$ with first marginal $\mu$ and second marginal $\nu$.  As a consequence  of Kantorovitch duality (see for instance \cite{Dud89}, chapter~11.8 or
 \cite{Vil03},
chapter~2), an equivalent definition of $W_2$ is
\begin{equation}\label{defwass2}
W^2_2(\mu,\nu)= \sup_{\phi \in \Xi }J(\phi):=\int_{\cK}\phi \mathrm d\mu + \int_{\cK} \phi^* \mathrm d\nu, \end{equation}
where $\Xi$ is the set of continuous functions $\phi \in  L^1_\mu(\cK, \mathbb{R})$ such that  $\phi(x)+\phi^*(y) \leq \|x-y \|^2,\ \mu\otimes\nu$-as with, for some arbitrarily chosen and fixed $x^* \in \cK$,
\[ \phi^*(x)=\inf_{y \in \cK}\|x-y\|^2-\phi(y), \mathrm{ \ } \phi=\left(\phi^*\right)^* \mathrm{\ and \ } \phi(x^*)=0.
\]

The supports of $\mu$ and $\nu$  are compact, so  any function $\phi$ in  $ \Xi$
is $2\|\cK\|$-Lipschitz, where $\|\cK\|$ is the diameter of $\cK$.  Thus Arzela-Ascoli's theorem implies that $\left(\Xi,\|\|_{\infty}\right)$ is relatively compact. Consequently  the  supremum  in formula  (\ref{defwass2}) is achieved;  we denote by $\Phi(\mu,\nu)$ the subset of $\Xi$ that maximizes $J(\phi,\phi^*)$. Its elements are called {\em Kantorovitch potentials} from $\mu$ to $\nu$.

\subsection{A differential game in Wasserstein space}
In the Wasserstein space of probability measures, the associated differential game we  consider is described as follows. For technical reasons -  the unicity of Kantorovich potentials -, the action sets of the players are not going to be $\Delta(X)$ and $\Delta(Z)$, but subsets of measures with a positive density lower bounded away from zero.

More precisely, let $\delta >0$ be some fixed parameter. Then there exist  \cite{Dud89} some  sets $\Delta_{\delta}(X) \subset \Delta(X)$ and $\Delta_{\delta}(Z)\subset \Delta(Z)$  such that
\begin{itemize}
\item[i)] $\Delta_{\delta}(X)$ is a convex and compact set;
\item[ii)] For every $\mu \in \Delta(X)$, there exists $\mu_\delta \in \Delta_{{\delta}}(X)$ such that $W_2(\mu,\mu_\delta) \leq \delta$;
\item[iii)] There exists $\underline{\delta}>0$ such that every $\mu_\delta \in \Delta_{{\delta}}(X)$ has a positive density lower-bounded by  $\underline{\delta}$.
\item[iv)] With a slight abuse of notations,  $\Delta_\delta(X\times Z)$ is the closed convex hull of the set of product measures $\mu\otimes \nu$ where $\mu \in \Delta_\delta(X)$ and $\nu \in \Delta_\delta(Z)$.
\end{itemize}

\medskip
We consider the differential game defined on $[s_0,1]$ with $s_0 \in (0,1)$ where
\begin{description}
\item{\textbf{Controls of players}}  are measurable maps $\bx$ and $\bz$ from $[s_0,1]$ to $\Delta_{{\delta}}(X)$ or $\Delta_{{\delta}}(Z)$; they are elements of $\cX_\delta(s_0)$ and $\cZ_\delta(s_0)$
\item{\textbf{NADC strategies}} are mappings from $\cZ_\delta(s_0)$ into $\cX_\delta(s_0)$ (for the player 1) satisfying the property given in Definition \ref{DF:NADCstrat} in Section \ref{SE:NADCStrat}. They are elements of $ \mathcal{A}_{NADC,\delta}(s_0)$ and similarly of $ \mathcal{B}_{NADC,\delta}(s_0)$ for player 2.

Because they are non-anticipative, the game can be written in normal form.
\item{\textbf{The dynamics}} are given  by  the integral formula:
\[ \theta(s) = \frac{s_0}{s} \theta_0+\frac{1}{s}\int_{s_0}^s \bx(t)\otimes\bz(t)dt, \ \theta(s_0)=\theta_0 \in \DD(X\times Z)\]
\item{\textbf{The terminal loss}} is $W^2_2\Big(\theta^{s_0,g_0,\alpha,\beta}(1),\widetilde{E}\Big)$  where $W_2^2(\cdot,\widetilde{E})$ is the square Wasserstein distance to a closed set $\widetilde{E}$.
\item{\textbf{The upper and lower values}} are the mapping defined by, for the upper value, \[V^+(s_0,\theta_{0})=\inf_{\alpha}\sup_{\beta}W_2^2(\theta^{s_0,\theta_0,\alpha,\beta}(1),\widetilde{E})=\inf_{\alpha}\sup_{\bz}W_2^2(\theta^{s_0,\theta_0,\alpha(\bz),\bz}(1),\widetilde{E}).\]
and, for the lower value, by
\[V^-(s_0,\theta_{0})=\sup_{\beta}\inf_{\alpha}W_2^2(\theta^{s_0,\theta_0,\alpha,\beta}(1),\widetilde{E})=\sup_{\beta}\inf_{\bx}W_2^2(\theta^{s_0,\theta_0,\bx,\beta(\bx)}(1),\widetilde{E}).\]
\item{\textbf{Name of the game.}} We shall denote this game as $\Gamma_{\widetilde{E},\delta}(s_0,\theta_0)$.
\end{description}

\subsection{Sub and super-solutions of Hamilton-Jacobi-Bellman equation in Wasserstein space. Comparison principle}

In this section we define  sub and super solutions of Hamilton-Jacobi-Bellman equation (HJB equation for short) in Wasserstein space and we obtain a comparison principle for this HJB equation. Let us first  define an adapted  concept of sub- and super-differential.

\medskip

\begin{definition} Let $\omega : [s_0,1] \times \DD(X \times Z) \to \R$ be a function and let $(\bar{t},\bar{\mu}) \in (s_0,1) \times \DD(X\times Z)$. We say that the pair $(p_{\bar{t}},\phi_{\bar{\mu}}) \in \R \times \Xi$ belongs to the super-differential $D^+\big(\omega(\bar{t},\bar{\mu})\big)$ to $\omega$ at $(\bar{t},\bar{\mu})$ if
\[ \limsup_{ \alpha \to 0, t \to \bar{t}} \sup_{\mu \in \Delta_\delta(X\times Z)}\frac{\omega(t, (1-\alpha)\bar{\mu}+\alpha\mu)-\omega(\bar{t},\bar{\mu})-p_{\bar{t}}(t-\bar{t})- \alpha \int \phi_{\bar{\mu}} \d (\mu-\bar{\mu})}{\alpha+|t-\bar{t}|} \leq 0.
\]
A pair $(p_{\bar{t}},\phi_{\bar{\mu}}) \in \R\times \Xi$ belongs to the sub-differential  $D^-\big(\omega(\bar{t},\bar{\mu})\big)$ to $\omega$ at $(\bar{t},\bar{\mu})$ if $(-p_{\bar{t}},-\phi_{\bar{\mu}}) \in \R\times \Xi$ belongs to the super-differential  $D^+\big(-\omega(\bar{t},\bar{\mu})\big)$ to $-\omega$ at $(\bar{t},\bar{\mu})$.
\end{definition}
\medskip

The supremum over $\mu \in \Delta_\delta(X\times Z)$ in this definition is the counterpart of the classical uniform convergence with respect to all possible directions in Euclidean spaces. 
Given an Hamiltonian $\cH\colon [s_0,1]\times \DD(X\times Z) \times \Xi \to \R $, we consider the associated  HJB equation: \begin{equation}\label{EQ:HJB}\omega_t+\cH(t,\mu,D\omega)=0\end{equation}
Its solutions are defined as follows.
\medskip

\begin{definition}
A sub-solution of the HJB Equation (\ref{EQ:HJB}) is an upper-semicontinuous  map $\omega : [s_0,1] \times \DD(X \times Z) \to \R$ such that  for any $(t,{\mu}) \in (s_0,1) \times \DD(X\times Z)$ and any $(p_{{t}},\phi_{{\mu}}) \in D^+\big(\omega({t},{\mu})\big)$ we have
\[ p_{{t}}+\cH({t},\mu,\phi_{{\mu}}) \geq 0.
\]
Super-solution are defined similarly.
\end{definition}
\medskip

We impose some regularity assumptions on  $\cH$ so that a comparison principle can be derived.
\medskip

\begin{hypo}\label{AS:Hamil} Assumptions on the regularity of $\cH$:
\begin{enumerate}
\item[i)] For any $\mu,\nu \in \DD(X \times Z)$, if $\phi$ is the Kantorovitch potential from $\mu$ to $\nu$ then
\[t\cH(t,\nu,-\phi^\star) - s\cH(s,\mu,\phi) \geq W_2^2(\mu,\nu)
\]
\item[ii)] $\cH$ is positively homogenous in $\phi$.
\end{enumerate}
\end{hypo}
\medskip

Using this assumption, we can derive the following comparison principle.
\begin{theorem}\label{TM:CompPrinc}\textbf{Comparison principle}

If $\omega_1$ and $\omega_2$ are respectively Lipschitz sub- and super-solution of Equation (\ref{EQ:HJB}) and    Assumption \ref{AS:Hamil} is satisfied then
\[ \inf_{[s_0,1]\times \DD(X\times Z)} (\omega_2-\omega_1)= \inf_{ \DD(X \times Z)} \omega_2(1,\cdot)-\omega_1(1,\cdot)=:A
\]
\end{theorem}
\medskip

\textbf{Proof.}
Let $k>0$ be a Lipschitz constant of  $\omega_1$ and $\omega_2$. Assume $A=0$ and let $(t_0,\mu_0)$ such that
\[ -\xi := \inf_{[s_0,1]\times \DD(X\times Z)} (\omega_2-\omega_1) <0 \mathrm{\ and\ } (\omega_2-\omega_1)(t_0,\mu_0) < - \xi/2
\]
and choose $\eta>0$ and $\gamma>0$ such that
\[\xi >2(\eta +k^2\gamma)\ \mathrm{and} \ \eta > \frac{2 k^2\gamma}{s_0}.\]

Let $\Phi$ on $\big([s_0,1]\times \DD(X\times Z)\big)^2$ be defined by
\[ \Phi(s,\mu,t,\nu)= - \omega_1(s,\mu)+\omega_2(t,\nu)+\frac{1}{\gamma}\Big(W_2^2(\mu,\nu)+(t-s)^2\Big)-\eta s
\]
and let $(\bar{s},\bar{\mu},\bar{t},\bar{\nu})$  be any of its minimizers.  The fact that $\Phi(\bar{s},\bar{\mu},\bar{t},\bar{\nu})\leq \Phi(\bar{s},\bar{\mu},\bar{s},\bar{\mu})$ implies that
\[ \frac{1}{\gamma}\Big(W_2^2(\bar{\mu},\bar{\nu}) + |\bar{s}-\bar{t}|^2\Big) \leq \omega_2(\bar{s},\bar{\mu})-\omega_2(\bar{t},\bar{\nu}) \leq k \sqrt{W_2^2(\bar{\mu},\bar{\nu}) + |\bar{s}-\bar{t}|^2 }. \]
where the last inequality follows from the fact that $\omega_2$ is $k$-Lipschitz. As a consequence, we immediately obtain that
\[
 W_2^2(\bar{\mu},\bar{\nu}) + |\bar{s}-\bar{t}|^2 \leq k^2 \gamma^2
\]

Assume that $\bar{s},\bar{t} \in (s_0,1)$ and let $\phi$ be the Kantorovitch potential from $\bar{\mu}$ to $\bar{\nu}$. Let $\mu \in \DD(X\times Z)$ be some fixed measure and, for every $\alpha \in (0,1]$, let $\phi_\alpha \in \Xi$ be the Kantorovitch potential from $(1-\alpha)\bar{\mu}+\alpha\mu$ to $\bar{\nu}$.   We recall that the mapping that associates  to a pair $(\mu,\nu) \in \DD(X\times Z)^2$ the set of Kantorovitch potentials is single-valued  because $\mu$ and $\nu$ have a density bounded away from zero. This mapping is also uniformly continuous as $\DD(X\times Z)$ is compact. We denote by $\Omega(\cdot)$ its modulus of continuity.

\bigskip

Since $\Phi(\bar{s},\bar{\mu},\bar{t},\bar{\nu})\leq \Phi(s,(1-\alpha)\bar{\mu}+\alpha\mu,\bar{t},\bar{\nu})$, we obtain that
 \[
 \omega_1(s,(1-\alpha)\bar{\mu}+\alpha\mu) \leq \omega_1(\bar{s},\bar{\mu})+ \frac{1}{\gamma}\Big(W_2^2\big((1-\alpha)\bar{\mu}+\alpha\mu,\bar{\nu}\big)-W_2^2(
 \bar{\mu},\bar{\nu})+(s-\bar{t})^2-(\bar{s}-\bar{t})^2\Big)+\eta(\bar{s}-s)
 \]
In particular, using the definition of $W_2$ in terms of Kantorovich potentials, \begin{eqnarray*}
W_2^2\big((1-\alpha)\bar{\mu}+\alpha\mu,\bar{\nu}\big)-W_2^2(
 \bar{\mu},\bar{\nu}) & \leq  \int \phi_\alpha \d \big\{(1-\alpha)\bar{\mu}+\alpha\mu\big\} + \int \phi_\alpha^\star \d \bar{\nu} -  \int \phi_\alpha \d \bar{\mu} - \int \phi_\alpha^\star \d \bar{\nu}\\
 &= \alpha \int \phi_\alpha \d (\mu - \bar{\mu})\\
 &= \alpha \int \phi \d (\mu - \bar{\mu})+ \alpha \int  (\phi_\alpha - \phi) \d (\mu - \bar{\mu})\\
 & \leq \alpha \int \phi \d (\mu - \bar{\mu}) + \alpha K_1 \|\phi_\alpha-\phi\|_\infty \\
 & \leq \alpha \int \phi \d (\mu - \bar{\mu}) + \alpha K_1\Omega(K_2\sqrt{\alpha})
 \end{eqnarray*}
 where   $K_1,K_2>0$ are constants depending on $X$ and $Z$.

\bigskip

Using the simple fact that \[(s-\bar{t})^2-(\bar{s}-\bar{t})^2=(s-\bar{s})(s+\bar{s}-2\bar{t})=2(s-\bar{s})(\bar{s}-\bar{t})+(s-\bar{s})^2\]
we therefore obtain that
\[ \left(\frac{2}{\gamma}(\bar{s}-\bar{t})-\eta,\frac{1}{\gamma}\phi\right) \in D^+\omega_1(\bar{s},\bar{\mu}) .
 \]
A similar proof give the dual result, i.e., $ \left(\frac{2}{\gamma}(\bar{s}-\bar{t}),-\frac{1}{\gamma}\phi^\star\right) \in D^-\omega_2(\bar{t},\bar{\nu})$.
\bigskip

Since  $\omega_1$ and $\omega_2$ are respectively sub and super-solution, we  therefore deduce that
\[ \frac{2}{\gamma}(\bar{s}-\bar{t})-\eta+\cH\left(\bar{s},\bar{\mu},\frac{1}{\gamma}\phi\right) \geq 0\ \mathrm{\ and\ }\ \frac{2}{\gamma}(\bar{s}-\bar{t})+\cH\left(\bar{t},\bar{\nu},-\frac{1}{\gamma}\phi^\star\right) \leq 0.
\]
The homogeneity and the regularity of $\cH$ yield that
\[ W_2^2(\bar{\mu},\bar{\nu}) \leq \bar{t} \cH(\bar{t},\bar{\nu},-\phi^\star)-\bar{s}\cH(\bar{s},\bar{\mu},\phi) \leq 2(\bar{s}-\bar{t})^2-\eta\bar{s}\gamma
\]

thus, dividing by $\gamma \bar{s}$,
\[ \eta \leq \frac{1}{\bar{s}\gamma}(-W_2^2(\bar{\mu},\bar{\nu})+2(\bar{s}-\bar{t})^2) \leq 2k^2\frac{\gamma}{s_0}\]
which is in contradiction with the choice of $\eta$.

\bigskip

It remains to check that $\bar{s}$ and $\bar{t}$ cannot be equal to $s_0$ or $1$ and this can be done exactly as in \cite{CarQui08}. For the sake of completeness, we provide the proof nonetheless.

Assume that $\bar{s}=1$ (the case $\bar{t}=1$ is identical). By definition of $\bar{\mu}, \bar{s}, \bar{\nu}, \bar{t}$ and $\mu_0,t_0$
\[ \Phi(\bar{s},\bar{\mu},\bar{t},\bar{\nu})\leq \Phi({t}_0,{\mu_0},{t_0},{\mu_0})=\omega_2(t_0,\mu_0)-\omega_1(t_0,\mu_0)-\eta t_0 \leq -\xi/2.
\]
Since $\bar{s}=1$ and $\omega_2$ is $k$-Lipschitz, we deduce  that
\[ \omega_2(1,\bar{\nu})-k\sqrt{W_2^2(\bar{\mu},\bar{\nu}) + |\bar{s}-\bar{t}|^2}-\omega_1(1,\bar{\mu})+\frac{1}{\gamma} \Big( W_2^2(\bar{\mu},\bar{\nu}) + |\bar{s}-\bar{t}|^2\Big)-\eta \leq - \xi/2
\]
The assumption that $A=\inf \omega_2(1,\cdot)-\omega_1(1,\cdot)=0$ yields that
\[  -k\sqrt{W_2^2(\bar{\mu},\bar{\nu}) + |\bar{s}-\bar{t}|^2}+\frac{1}{\gamma} \Big( W_2^2(\bar{\mu},\bar{\nu}) + |\bar{s}-\bar{t}|^2\Big)-\eta \leq - \xi/2
\]
which is impossible given the choice of $\xi$. If $\bar{s}=s_0$, then we conclude using the fact that a sub- or super-solution on $(s_0,1]$ is a sub- or super-solution on $[s_0,1]$.
$\hfill \Box$

\bigskip

The proof of Theorem \ref{TM:CompPrinc} indicates that condition \textit{i)} of Assumption \ref{AS:Hamil} could actually be replaced by the following weaker version. There exists some $k'>0$ such that \[\label{qq} t\cH(t,\nu,-\phi^\star)-s\cH(s,\mu,\phi) \geq -k' W_2^2(\mu,\nu), \ \forall \mu,\nu \in \Delta_\delta(X\times Z).\]\subsection{Existence of a value in the differential game}

As in Euclidean space, we can derive some regularity of the upper and the lower  value functions.
\medskip

 \begin{proposition}\label{PR:RegV} Both the upper and the lower value functions $V^+(\cdot,\cdot)$ and $V^-(\cdot,\cdot)$ are Lipschitz on $(s_0,1]\times \DD(X\times Z)$,  for every $s_0 >0$, with a Lipschitz-constant independent of  $\delta$ and $s_0$. They can be uniquely extended to $[0,1]$ to mappings such that $V^{\pm}(0,\cdot)$ are constant.

Moreover, the upper value function $V^+$ is a sub-solution to (\ref{EQ:HJB}) with \[\cH=\cH^+(t,\mu,\phi):=\inf_{\bx}\sup_{\bz} \frac{1}{t} \int \phi \d (\bx\otimes\bz-\mu).\]
 The lower value function $V^-$ is a super-solution to (\ref{EQ:HJB}) with \[\cH=\cH^-(t,\mu,\phi):=\sup_{\bz}\inf_{\bx} \frac{1}{t} \int \phi \d (\bx\otimes\bz-\mu).\]
 \end{proposition}
 \medskip

\textbf{Proof.}
 Using the exact same proofs than in Appendix \ref{SE:Appendix}, we can show that  $V^+$ and $V^-$ are Lipschitz and they satisfy the dynamic programming principle. It only remains to show that it implies that they are sub- and super-solution of (\ref{EQ:HJB}), which was a well known fact in Euclidean space (see Lemma \ref{DPP} and Theorem \ref{egalS}).

  \medskip

The dynamic programming principle and the integral form of $\theta(s)$ imply that  for $0 < t_0 < t_0+h < 1$
 \begin{equation}\label{dp} V^+(t_0,\theta_0)=\inf_\alpha \sup_{\mathbf{z}} V^+\left(t_0+h,(1-\varepsilon_{t_0,h})\theta_0+\varepsilon_{t_0,h}\lambda_{\alpha,\mathbf{z},t_0,h}\right) \end{equation} where $\displaystyle\lambda_{\alpha,\bz,t_0,h}:=\frac{1}{h}\int_{t_0}^{t_0+h} \alpha (\bz) (s)\otimes\bz(s) \d s$  and $ \displaystyle \varepsilon_{t_0,h}:=\frac{h}{t_0+h}$.

Take $ \alpha$ a NADC strategy.  For $h$ small enough, the strategy  $\alpha $ is constant on $(t_0, t_0+h)$. So  there exists some $\bx _ \alpha \in \Delta_\delta(X)$ such for any $ \bz ( \cdot) $ we have $ \alpha  (\bz) (s) =  \bx _ \alpha$ for any $s \in (t_0, t_0+h)$.
 
  Let $(p_t,\phi_\mu) \in D^+V^+(t_0,\mu_0)$, then the definition of sub differential implies that for any $\bz (\cdot) $ 
 \begin{align*} & \frac{ V^+(t_0+h,(1-\varepsilon_{t_0,h})\theta_0+\varepsilon_{t_0,h}\lambda_{\alpha,\mathbf{z},t_0,h})-V^+(t_0,\theta_0)}{h} \\ & \hspace{1cm} \leq 
 p_t+ \frac{1}{h} \varepsilon_{t_0,h}\int \phi_\mu \d (\lambda_{\alpha,\mathbf{z},t_0,h}-\theta_0) 
 +  \frac{1}{h}  (h+\varepsilon_{t_0,h})o(h+\varepsilon_{t_0,h})
 \end{align*} where $ o(h+\varepsilon_{t_0,h})  \to 0 $ as  $h \to 0 ^+$ uniformly with respect to  $ \bz(\cdot)$ and $\alpha$.   Since $\alpha $ is constant on $(t_0, t_0+h)$, we get
 \begin{align*} & \frac{ V^+(t_0+h,(1-\varepsilon_{t_0,h})\theta_0+\varepsilon_{t_0,h}\lambda_{\alpha,\mathbf{z},t_0,h})-V^+(t_0,\theta_0)}{h} \\ & \hspace{1cm}  \leq
 p_t+ \varepsilon_{t_0,h}\int _{X\times Z}\phi_\mu \d (\frac{1}{h}\int_{t_0}^{t_0+h} \bx _ \alpha \otimes\bz(s) \d s-\theta_0) 
 +  \frac{1}{h}  (h+\varepsilon_{t_0,h})o(h+\varepsilon_{t_0,h}) \\
 &\hspace{1cm}  \leq  p_t+ \varepsilon_{t_0,h}\sup_{\bz  \in   \Delta_\delta(Z)}\int _{X\times Z}\phi_\mu \d ( \bx _ \alpha \otimes\bz-\theta_0) 
 +  \frac{1}{h}  (h+\varepsilon_{t_0,h})o(h+\varepsilon_{t_0,h})
 \end{align*} 
 Observe that  any constant control $ \bx$  can generate a NADC strategy $ \alpha$  such that $\alpha $ is constant equal to $\bx$ on $(t_0, t_0+h)$. So passing to the supremum over $\bz(\cdot)$ and to the infimum over $ \alpha$, we obtain in view of \eqref{dp} \begin{eqnarray*} 0
  \leq  p_t+ \varepsilon_{t_0,h_n} \inf  _ {\bx  \in   \Delta_\delta(Z)} 
 \sup_ {\bz  \in   \Delta_\delta(Z)} \int _{X\times Z}\phi_\mu \d (\bx \otimes \bz -\theta_0) 
 +  \frac{  h_n+\varepsilon_{t_0,h_n}}{h_n}o(h_n+\varepsilon_{t_0,h_n}) .
 \end{eqnarray*}
Letting $h \to 0$, this gives  \[ 0 \leq  p_t+ \inf_{\bx \in \Delta_\delta(X)}\sup_{\bz \in \Delta_\delta(Z)} \frac{1}{t_0} \int \phi_\mu \d (\bx\otimes\bz-\theta_0)\]  thus $V^+$ is a sub-solution of  (\ref{EQ:HJB}) with respect to the Hamiltonian defined by $\cH^+(t,\mu,\phi):=\inf_{\bx}\sup_{\bz} \frac{1}{t} \int \phi \d (\bx\otimes\bz-\mu)$.

\bigskip

 The proof that $V^-$ is a super-solution is similar and  is omitted. $\hfill \Box$

 \bigskip

Before being able to state the existence of the value by using the comparison principle stated in Theorem \ref{TM:CompPrinc}, we need to prove that $\cH^+$ and $\cH^-$ satisfy Assumption~\ref{AS:Hamil}.
\begin{proposition}
  $\cH^+=\cH^-$ and they satisfy Assumption \ref{AS:Hamil}.
 \end{proposition}
 \textbf{Proof.}  Assume that $\phi$ is the Kantorovitch potential from $\mu$ to $\nu$, then  \begin{eqnarray*} t\cH^+(t,\mu,\phi)&=&\inf_{\bx}\sup_{\bz} \int \phi \d ( \bx\otimes\bz-\mu)\\&=&-\int \phi\d\mu-\int\phi^\star\d\nu+\inf_{\bx}\sup_{\bz} \int\phi^\star\d\nu+\int\phi\d\
bx\otimes\bz\\
 &\leq& -d(\mu,\nu)^2+\inf_{\bx}\sup_{\bz} \int (-\phi^*)\d (\bx\otimes\bz-\nu)\\&=&-d(\mu,\nu)^2+s\cH^+(s,\mu,-\phi^\star),
  \end{eqnarray*}
where the inequality is a consequence of the fact that  $\phi \leq - \phi^*$. So $\cH^+$ satisfies Assumption \ref{AS:Hamil}.

\medskip

Since any mapping $\phi \in \Xi$ is continuous and $X,Z$ are compact sets in Euclidean space, Sion \cite{Sio58} minmax theorem implies that $\cH^+=\cH^-$. $\hfill \Box$

  \bigskip

  \begin{corollary}
The game   $\Gamma_{\widetilde{E},\delta}(s_0,\theta_0)$ has a value $V$. Moreover,  $V(0,\cdot)$ is constant.
 \end{corollary}
\textbf{Proof.}   The existence is a direct consequence of the comparison principle stated in Theorem \ref{TM:CompPrinc} and the fact that $V$ is constant is due to the regularity property of $V^\pm$ stated in Proposition \ref{PR:RegV}. $\hfill \Box$
\medskip

We denote by $V_{\widetilde{E},\delta}$ the value of the constant mapping $V(0,\cdot)$.

  \subsection{From  differential game  to repeated game and weak approachability}
 As in the Euclidian case, we can now formulate the general approachability theorem:
 \medskip

 \begin{theorem}
 Any closed set $\widetilde{E} \subset \Delta(X\times Z)$ is either weakly-approachable or weakly-excludable.

 More precisely, $\widetilde{E}$ is weakly approachable if and only if $\sup_{\varepsilon}\lim\inf_{\delta} V_{\widetilde{E}^\varepsilon,\delta}=0$.
 \end{theorem}\medskip

\textbf{Proof.}
 Assume first that  $\sup_{\varepsilon}\lim\inf_{\delta} V_{\widetilde{E}^\varepsilon,\delta}=0$ and let $\varepsilon$ be fixed and such that $\lim\inf_{\delta} V_{\widetilde{E}^\varepsilon,\delta}=0$.

 For any $\varepsilon'>0$, let  $\delta>0$ be such that $V_{E^\varepsilon,\delta}<\varepsilon'$, $s_0 \leq \varepsilon'/L$ where $L$ is the Lipschitz constant of $V$ and $\theta_0\in\DD(X\times Z)$ chosen arbitrarily.

 \bigskip

Given  any  ${\varepsilon'}$-optimal strategy $\alpha \in \mathcal{A}_{NADC}(s_0,\theta_0)$ in  the game $\Gamma_{\widetilde{E}^\varepsilon,\delta}(s_0,\theta_0)$,  we are going to construct a strategy $\sigma_{\alpha,n}$ of  player 1 in the repeated game.  This is done  almost  exactly as in Section \ref{SE:FromDiffToRep}. The only difficulty is that  player 2 can choose at stage $n \in \N$ some action $\bz_n \in \Delta(Z)$ that does not belong to $ \DD(Z)$. In that case, we approximate $\bz_n$ by $\bz_n^{\delta} \in \DD(Z)$ such that $W_2(\bz_n^\delta,\bz_n) \leq \delta\leq {\varepsilon}$ and the strategy $\sigma_{\alpha,n}$ is defined with respect to $\alpha$ and the sequence $\{\bz^\delta_n\}$.

\bigskip

Denoting $\theta_n^\delta=\bx_n \otimes \bz_n^\delta$, the strategy $\sigma_{\alpha,n}$  we have constructed is such that for any $n \in \N$ large enough (see the proof of Theorem \ref{TH:Vieille}),
\[ W^2_2(\overline{\theta}_n^\delta,\widetilde{E}) \leq V_{\widetilde{E}^\varepsilon,\delta} + 4\varepsilon' \leq 5 \varepsilon'.
\]
This entails the weak-approachability of $\widetilde{E}$ since (for $n$ large enough)
\[ W_2(\overline{\theta}_n,\widetilde{E})\leq W_2(\overline{\theta}_n,\widetilde{E}^\varepsilon)+\varepsilon  \leq W_2(\overline{\theta}_n^\delta,\widetilde{E}^\varepsilon)+\varepsilon + \delta \leq \sqrt{5 \varepsilon'}+\varepsilon+\delta.
\]
Letting $\varepsilon'$ and $\delta$  to zero entails the result.

\bigskip

Reciprocally, assume that  $\sup_{{\varepsilon}}\lim\inf_{{\delta}} V_{E^{\varepsilon},{\delta}}=\eta >0$. Thus for some $\varepsilon>0$ and every $\delta$ small enough, $V_{E^{\varepsilon},{\delta}} \geq \eta/2$. This implies that player 2 can weakly approach $\widetilde{E}^{\varepsilon}$ (as before, if  player 1 chooses $\bx_n \not \in \DD(X)$, player 2 can respond as if he played $\bx_n^{\delta}\in \DD(X)$).

\bigskip

Therefore any closed set is either weakly-approachable  of weakly-excludable.
$\hfill \Box$

\bigskip

This result has an important corollary in repeated game with partial monitoring:
\medskip

\begin{corollary}
Any closed set $E \subset \R^d$ such that there exists a some closed set $\widetilde{E} \subset \Delta(X\times Z)$ satisfying \[ E=\{\mathbb{E}_{q}[\mathbf{p}] , q \in \widetilde{E} \} \subset \R^d
\]
 is either weakly-approachable or weakly-excludable. If there exists no such set $\widetilde{E}$, then it is possible that $E$ is neither weakly-approachable nor weakly-excludable (and, furthermore, $E$ can even be a convex and compact subset of $\R^d$).
\end{corollary}

\section{Appendix}\label{SE:Appendix}

In this section, we recall some well known facts on differential games. We also  prove  that the values of a rather general differential game with nonanticipative strategies with delay coincide with the values defined with NADC strategies. This is a new result with is valuable by itself independently on the problem of weak approachability.

Let us consider the dynamics
 \begin{equation} \label{eqgen} g'(s) = f(s, g(s), x(s), y(s)) , \; s \in  [a,1] \end{equation}  where $a  <1$ is fixed and
 $ f : [a,1] \times \R ^d\times   X \times Y  \mapsto \R ^ d$ is a Lipschitz continuous bounded map. Let $ \ell : \R ^d \mapsto \R$ be a Lipschitz continuous function. For any $ s \in [a, 1]  $ $ g_ 0 \in \R ^d $ and
  $(\bx, \by) \in \cX (s_0) \times \cY (s_0) $  we denote by $ s \mapsto g^{s_0,g_0,\bx, \by}(s) $ the unique solution to (\ref{eqgen}). Similarly for any pair $(\alpha, \beta) \in 
  \mathcal{A} (s_0)  \times
  \mathcal{B} (s_0) $
 we set   $  g^{s_0,g_0,\alpha, \beta } = g^{s_0,g_0,\alpha, \beta} $ where  $(\bx, \by) $ is associated with  $(\alpha, \beta) $ by Lemma \ref{rep}.

For every $s_0 \in (0,1]$ and $g_0 \in \R ^d$,  one can define
the upper-value of the game
\[V^+(s_0,g_{0}):=\inf_{\alpha \in  \mathcal{A} (s_0) }\sup_{\beta \in  \mathcal{B} (s_0) }\ell(g^{s_0,g_0,\alpha,\beta}(1)).\]
and  the lower-value
\[V^-(s_0,g_{0})=\sup_{\beta \in  \mathcal{B} (s_0) }\inf_{\alpha \in  \mathcal{A} (s_0) }\ell(g^{s_0,g_0,\alpha,\beta}(1)).\]
Lemma \ref{rep}  implies  the following:
\[V^+(s_0,g_{0}) =\inf_{\alpha \in  \mathcal{A} (s_0) }\sup_{\by}\ell(g^{s_0,g_0,\alpha(\by),\by}(1)), \; V^-(s_0,g_{0})=\sup_{\beta \in  \mathcal{B} (s_0) }\inf_{\bx}\ell(g^{s_0,g_0,\bx,\beta(\bx)}(1)).\] Due to the Lipschitz continuity of the cost function $\ell$ and the regularity of the dynamics~\eqref{eqgen}, one can easily obtain from the above relation that $V^+$ and $V^-$ coincide with the values defined in subsection \ref{basic}.

Under the above assumption it is well known that $ V ^+ $ and  $ V ^- $ are Lipschitz continuous and they are respectively viscosity solutions of Hamilton Jacobi Isaacs Equations       (cf for instance  \cite{ES,  BCD} see also \cite{BCQ} and the references therein): $V ^+ $ solves
\begin{equation}\label{HJ+} \left\{\begin{array}{ll} \frac{\partial V}{\partial s}(s,g)+H^{+}\left(s,g,\frac{\partial V}{\partial g}(s,g)\right) =0 & \mathrm{for all }\ (s,g) \in [a,1]\times \R^d\\
V(1,g)=\ell(g) & \mathrm{ for all }\ g \in \R ^d\end{array}\right.
\end{equation}
while $V^-$ is a solution to
\begin{equation}\label{HJ-}
 \left\{\begin{array}{ll} \frac{\partial V}{\partial s}(s,g)+H^{-}\left(s,g,\frac{\partial V}{\partial g}(s,g)\right) =0 & \mathrm{ for all }\ (s,g) \in [a,1]\times \R^d \\
V(1,g)=\ell(g) & \mathrm{for all }\ g \in \R ^d \end{array}\right.
\end{equation}
where the Hamiltonians are defined respectively by
\[H^-(s,g,p):=\sup_{y\in Y} \inf_{x \in X}p.f(s,g,x,y)  \ \quad \mathrm{and}\quad H^+(s,g,p):= \inf_{x \in X}\sup_{y\in Y}p.f(s,g,x,y)  .\]
Moreover the PDE (\ref{HJ+}) has a unique bounded uniformly continuous viscosity solution \cite{CIL, BCD}, the same property holds for (\ref{HJ-}).

If furthermore we assume the following Isaacs equation
\begin{equation}\label{I}
H^-(s,g,p) =H^+(s,g,p) ,  \; \forall (s,g,p) ,
\end{equation} then the Hamilton Jacobi Equations (\ref{HJ+}) and  (\ref{HJ-}) are the same and consequently $ V ^+ = V ^- $ (the differential game has a  value).

One can also define the values using strategies NADC of Definition \ref{DF:NADCstrat} as follows

\[V_{NADC}^+(s_0,g_{0}):=\inf_{\alpha \in  \mathcal{A} _{NADC} ( s_0  ) }\sup_{\beta \in  \mathcal{B} _{NADC}( s_0  ) }\ell(g^{s_0,g_0,\alpha,\beta}(1)).\]

\[V_{NADC}^-(s_0,g_{0}):=\sup_{\beta \in  \mathcal{B}_{NADC} ( s_0  ) }\inf_{\alpha \in  \mathcal{A}_{NADC}( s_0 ) }\ell(g^{s_0,g_0,\alpha,\beta}(1)).\]
Once again observe that  Lemma \ref{rep} yields
\begin{align*}V_{NADC}^+(s_0,g_{0}) &=\inf_{\alpha \in  \mathcal{A}_{NADC} (s_0) }\sup_{\by}\ell(g^{s_0,g_0,\alpha(\by),\by}(1)), \\ V_{NADC}^-(s_0,g_{0}) &=\sup_{\beta \in  \mathcal{B}_{NADC} (s_0) }\inf_{\bx}\ell(g^{s_0,g_0,\bx,\beta(\bx)}(1)).\end{align*}

We now prove that under our assumptions the NADC strategies define the same values that values defined through NAD strategies.

\begin{theorem} \label{egalS}

$$  V_{NADC}^+ = V ^ + \mbox{and } V_{NADC}^ - = V ^-. $$
  $V_{NADC}^+ = V ^ +$ is the unique    bounded uniformly continuous  viscosity solution to (\ref{HJ+}) while  $V_{NADC}^- = V ^ -$ is the unique    bounded uniformly continuous  viscosity solution to (\ref{HJ-})

\end{theorem}
\textbf{Proof.}
We first show that the values  $V_{NADC}^\pm $ satisfy  the  following dynamic programming principle (which proof is postponed later on).
\begin{lemma}\label{DPP}  for all $  0<s_0 < s_1 \leq 1$ ,
\[ \left\{ \begin{array}{ll} i) &    V_{NADC}^-(s_0,g_0) =\sup_{\beta \in  \mathcal{B}_{NADC} ( s_0  ) }  \inf_{\bx} V_{NADC}^-(s_1,g^{s_0,g_0,\bx,\beta(\bx)}(s_1))
\\ ii ) &   V_{NADC}^+(s_0,g_0) = \inf_{\alpha \in  \mathcal{A}_{NADC} (s_0) } \sup_{\by} V_{NADC}^+(s_1,g^{s_0,g_0,\alpha (\by) ,\by)}(s_1))
\end{array}\right.
.\]
\end{lemma}
Once the dynamic programming  is obtained, it is standard \cite{BCD} to prove that $ V_{NADC}^+ $ and $  V_{NADC}^-$ are viscosity solution of (\ref{HJ+}) and (\ref{HJ-}) respectively.

As for $V^+$ and $V ^- $, the same arguments allow to obtain that   $ V_{NADC}^+ $ and $  V_{NADC}^-$  are Lipschitz continuous.

 $ V_{NADC}^+ $ and  $V^+$  are both viscosity solutions to
(\ref{HJ+})  which has  a unique bounded uniformly continuous viscosity solution by the comparison theorem of  \cite{CIL}).   So we deduce that $ V_{NADC}^+ = V ^ + $. Similarly we obtain $V_{NADC}^ - = V ^-$, which completes the proof of Theorem \ref{egalS}.
$\hfill \Box$

\medskip

It only remains to prove Lemma \ref{DPP}.

{\em Proof of Lemma \ref{DPP}} We only prove the part i) of the Lemma since part ii) is very similar.

Denote by $
W(s_0,s_1,g_0):=\sup_\beta\inf_{\bx} V^-(s_1,g^{s_0,g_0,\bx,\beta(\bx)}(s_1)) $

The proof goes in two steps, we first prove that $V_{NADC}^-(s_0,g_0) \leq W(s_0,s_1,g_0)$ and then that the converse inequality also holds.

\begin{itemize}
\item[i)]Let $\beta_0 \in \mathcal{B}_{NADC}(s_0)$, $\bx_0 \in \cX(s_0)$ and define $g(s)=g^{s_0,g_0,\bx_0,\beta_0(\bx_0)}(s)$. We construct a strategy $\beta_1$ by
\[ \forall \bx \in \cX(s_1), \ \beta_1(\bx)=\beta_0(\tilde{\bx}) \ \mathrm{ where } \ \tilde{\bx}=\left\{\begin{array}{ll} \bx_0 & \mathrm{on}\ [s_0,s_1]\\ \bx & \mathrm{on}\ [s_1,1]\end{array}\right.
\]
By construction, $\beta_1$ belongs to $\cB_{NADC}$ because $\beta_0$ does; $g^{s_0,g_0,\tilde{\bx},\beta_0(\tilde{\bx})}$ coincides with $g^{s_1,g(s_1),\bx,\beta_1(\bx)}$ on $[s_1,1]$ thus $\ell\Big(g^{s_1,g(s_1),\bx,\beta_1(\bx)}(1)\Big)= \ell\Big(g^{s_0,g_0,\tilde{\bx},\beta_0(\tilde{\bx})}(1)\Big)$ and
\begin{eqnarray*} \inf_{\bx \in \cX(s_1)} \ell\Big(g^{s_1,g(s_1),\bx,\beta_1(\bx)}(1)\Big)&=&\inf_{\bx \in \cX(s_0),\ \bx=\bx_0 \ \mathrm{on} \ [s_0,s_1]}\ell\Big(g^{s_0,g_0,\bx,\beta_0(\bx)}(1)\Big)  \\ &\geq& \inf_{\bx \in \cX(s_0)}\ell\Big(g^{s_0,g_0,\bx,\beta_0(\bx)}(1)\Big).
\end{eqnarray*}
The term on the left is smaller than $V^-(s_1,g(s_1))=V^-(s_1,g^{s_0,g_0,\bx_0,\beta_0(\bx_0)}(s_1))$; since $\bx_0$ and $\beta_0$ where taken arbitrarily, we deduce that
\[ \sup_{\beta_0 \in  \mathcal{B}_{NADC}}\inf_{\bx_0}V^-(s_1,g^{s_0,g_0,\bx_0,\beta_0(\bx_0)}(s_1)) \geq \sup_{\beta_0 \in  \mathcal{B}_{NADC}}\inf_{\bx} \ell\Big(g^{s_0,g_0,\bx,\beta_0(\bx)}(1)\Big)\]
which is the first inequality claimed.
\item[ii)] We now prove the converse inequality. Let $g_0$ and  $\varepsilon >0$ be fixed, $s_1^+$ be any rational number in $[s_1,s_1+\varepsilon)$. Consider $$R(g_0) := \{ g ^{s_0,g_0,\bx,\by } (s)  , \;
 (\bx ,\by)  \in  \cX(s_0) \times  \cY(s_0), \, s \in [s_0,1] \,\} $$ the attainable set starting from $(s_0,g_0)$ which is a precompact subset of
$  \R ^d$.  Given $g_1 \in R(g_0)$, let us define  $\beta_{g_1} \in \cB_{NADC}(s^+_1)$ such that
\[ \inf_{\bx \in \cX(s^+_1)} \ell\Big(g^{s^+_1,g_1,\bx,\beta_{g_1}(\bx)}(1)\Big) > V_{NADC}^-(s^+_1,g_1)-\varepsilon.
\]
Since $V_{NADC}^-$ is Lipschitz and $\overline{R(g_0)}$ is compact, we can assume that there exists only a finite number of different strategies $\beta_{g}$.

Given $\beta_0 \in \cB_{NADC}(s_0)$, we construct a new strategy $\beta\in\cB_{NADC}(s_0)$ by
\[ \forall \bx \in \cX(s_0),\ \beta(\bx)= \left\{\begin{array}{ll} \beta_0(\bx) & \mathrm{on}\ [s_0,s_1^+]\\ \beta_{g(s^+_1)}(\bx_{|_{[s_1^+,1]}}) & \mathrm{on}\ [s_1^+,1]\end{array}\right.\, ,
\]
with  $g(s^+_1)=g^{s_0,g_0,\bx,\beta(\bx)}(s^+_1)=g^{s_0,g_0,\bx,\beta_0(\bx)}(s^+_1)$. The fact that $\beta$  belongs to $\cB_{NADC}(s_0)$ comes from the fact that there exists only a finite number of $\beta_g$ and $s_1^+$ is rational (this ensures the existence of the integer $N$ required by  point 1)). Therefore, we obtain for any $\bx \in \cX(s_0)$,
\[ \ell\Big(s_0,g_0,\bx,\beta(\bx)\Big)=\ell\Big(s^+_1,g(s^+_1),\bx_{|_{[s^+_1,1]}},\beta_{g(s^+_1)}(\bx_{|_{[s^+_1,1]}})\Big)\geq V_{NADC}^-(s^+_1,g(s^+_1))-\varepsilon.
\]
Taking the infimum in $\bx \in \cX(s_0)$ and the supremum in $\beta_0 \in \cB_{NADC}(s_0)$ yield
\[ V_{NADC}^-(s_0,g_0) \geq \sup_{\beta_0}\inf_{\bx \in \cX(s_0)} V_{NADC}^-(s^+_1,g(s^+_1))-\varepsilon.
\]
Because $\|g(s_1+)-g(s_1)\| \leq \|f\|_\infty|s_1^+-s_1|\leq \|f\|_\infty\varepsilon$ and since $V_{NADC}^-$ is Lipschitz, we obtain that
\[V_{NADC}^-(s_0,g_0)\geq W(s_0,s_1,g_0)-(1+\|f\|_\infty)\Lip(V_{NADC}^-)\varepsilon-\varepsilon\, .\]
\end{itemize}
This completes the proof of Lemma \ref{DPP}.  $ \Box $

\bigskip

\underline{Remark}
The same idea of proof shows that Theorem \ref{egalS}  also holds true for values defined through the strategies defined below.
A fixed delay nonanticipative strategy for player I is a map $ \alpha : \cY(s_0) \mapsto \cX(s_0) $ such that there exists  $\tau >0 $   such that for any $ t \in  [s_0 ,1]  $ if $ \by _1 ( \cdot) )  $ and $ \by _2 ( \cdot) $ coincide almost surely on $ [s_0 , t]$ then the controls  $ \alpha (\by _1 ) ( \cdot) )  $ and $ \alpha (\by _2)  ( \cdot) $ coincide almost surely on $ [s_0 , \min\{t+ \tau , 1\} ]$.  We define in a similar way the set $ \mathcal{B}(s_0) $ of nonanticipative strategies with fixed delay  $ \beta $ for the other player.

 \medskip
 
 MPS14

\end{document}